\documentclass[12pt]{amsart}

\usepackage[utf8]{inputenc}
\usepackage[T2A]{fontenc}

\usepackage{graphicx}
\usepackage{amsmath}
\usepackage{amscd}
\usepackage{amssymb, amsmath, amsthm}
\usepackage[all]{xy}
\usepackage{enumitem}
\usepackage[margin=2.5cm]{geometry}

\newtheorem{theorem}[subsection]{Theorem}
\newtheorem{lemma}[subsection]{Lemma}

\newtheorem{definition}[subsection]{Definition}

\newtheorem{remark}[subsection]{Remark}

\makeatletter
\@addtoreset{subsection}{section}
\@addtoreset{equation}{section}
\@addtoreset{figure}{section}
\@addtoreset{table}{section}
\makeatother

\makeatletter
\newcommand\testshape{family=\f@family; series=\f@series; shape=\f@shape.}
\def\myemphInternal#1{\if n\f@shape%
\begingroup\itshape #1\endgroup\/%
\else\begingroup\bfseries #1\endgroup%
\fi}
\def\myemph{\futurelet\testchar\MaybeOptArgmyemph}
\def\MaybeOptArgmyemph{\ifx[\testchar \let\next\OptArgmyemph
             \else \let\next\NoOptArgmyemph \fi \next}
\def\OptArgmyemph[#1]#2{\index{#1}\myemphInternal{#2}}
\def\NoOptArgmyemph#1{\myemphInternal{#1}}
\makeatother

\newcommand\RRR{\mathbb{R}}
\newcommand\ZZZ{\mathbb{Z}}
\newcommand\PPP{\mathbf{P}}
\newcommand\FF{\mathcal{F}}

\newcommand\Rntriv{\underline{\RRR}^n}

\newcommand\id{\mathrm{id}}
\newcommand\diam{\mathrm{diam}}

\newcommand\Mman{M}
\newcommand\Nman{N}

\newcommand\xsol{\xi} 
\newcommand\axsol{\lambda} 

\newcommand\EE{\mathbf{E}}

\newcommand\Xmap{X}

\newcommand\Rinj[1]{R_{#1}}

\newcommand\mymid{\mid}

\title{Contractibility of manifolds by means of stochastic flows}

\author{Alexandra Antoniouk, Sergiy Maksymenko}
\email{antoniouk.a@gmail.com, maks@imath.kiev.ua}

\address{Institute of Mathematics of NAS of Ukraine, Te\-re\-shchenkivska st. 3, Kyiv, 01004 Ukraine}

\thanks{The authors would like to thank Vladimir Krouglov and Dmitry Bolotov for useful discussions.}



\begin{document}

\begin{abstract}
In the paper [Probab.\! Theory Relat.\! Fields, {\bf 100} (1994) 417--428] Xue-Mei Li has shown that the moment stability of an SDE is closely connected with the topology of the underlying manifold.
In particular, she gave sufficient condition on SDE on a manifold $M$ under which the fundamental group $\pi_1 M=0$. 
We prove that in fact under the similar conditions the manifold $M$ is contractible, that is all homotopy groups $\pi_k M$, $k\geq1$, vanish.
\end{abstract}

\maketitle

\section{Introduction}
The interplay of geometrical or topological structures of manifold and the properties of differential operations on it form a library of the most crucial results in analysis. 
For instance, 
\begin{enumerate}[leftmargin=*, label=$(\arabic*)$]
\item\label{enum::Morse_ineq}
if $\Mman$ is closed, then the number of (non-degenerate) critical points of index $i$ of a Morse function $f:\Mman\to\RRR$ bound the ranks of $i$-th homology groups of $\Mman$ (Morse inequalities);

\item\label{enum::deRham}
de Rham cohomologies $H^{*}_{\mathrm{DR}}(\Mman)$ of an orientable manifold $\Mman$ are isomorphic with singular real cohomologies $H^{*}(\Mman,\RRR)$, (de Rham theory);

\item\label{enum::Riemann_metric}
if $\Mman$ is Riemannian, then there is a lot of statements relating its Ricci and sectional curvatures with the topology of $\Mman$  and especially with the fundamental group $\pi_1\Mman$, (theorems by Cartan-Hadamard, Bonnet-Myers, Preissman, Byers, Bochner);

\item\label{enum::HopfPoincare}
for a vector field $F$ on $\Mman$ having only isolated zeros the alternating sum of indexes of those zeros equals the Euler characteristic of $\Mman$, (Poincar\'e-Hopf theorem);

\item\label{enum::TopEntropy}
topological entropy of smooth dynamical systems on $\Mman$ can also be computed via Lyapunov exponents (Margulis-Ruelle Inequality, Pesin entropy formula).
\end{enumerate}

The invention of the stochastic analysis since the milestone papers of Wiener and Ito gave rise to a question of natural extension of the above results to the new kind of differential objects connected with the stochastic theory to get from them an information about the geometry and topology of underlying manifold $\Mman$.
In 1962 Ito~\cite{Ito:1962} (see also \cite{Ito:1975Displ}) introduced a notion of a stochastic parallel transport, which generalizes the notion of parallel transport in differential geometry.
These ideas permitted further development of stochastic analysis on manifold developed in the papers of Eells and Elworthy \cite{EElw:1970}, Malliavin~\cite{Mal:1974}, Airault~\cite{Airault:1976}, Vauthier~\cite{Vauthier:1979}, Berthier and Gaveau \cite{BerthierrGaveau:1978}, and many others.
The principal problem which appears here is that the corresponding objects depend only 
\myemph{measurably} and not \myemph{continuously} on probabilistic parameters.
This essentially prevents the use of the topological logic of continuity and well developed homotopy invariants.

On the other hand, for a Riemannian manifold the Laplace operator $\Delta$ uniquely defines a Brownian motion on $\Mman$.
This allowed to prove analogues of results of type~\ref{enum::Riemann_metric} in terms of stochastic differential equation (SDE) on $\Mman$.

Another approach is based on extending results of type~\ref{enum::TopEntropy}.
Recall that a \myemph{maximal Lyapunov exponent} of a diffeomorphism $h:\Mman\to\Mman$ of a Riemannian manifold $\Mman$ at a point $x_0\in\Mman$ is defined by
\begin{equation}\label{equ:Lyapunov_exp}
\lambda(x_{0}) = \lim _{{n\to \infty }}{\frac {1}{n}}\sum_{{i=0}}^{{n-1}}\ln \|h'(x_{i})\| 
\end{equation}
where $\|h'(x)\|$ is the norm of the tangent linear map $T_x h: T_x \Mman \to T_{h(x)}\Mman$, and $x_i = h^i(x_0)$ is an $i$-th iteration of $x_0$ under $h$.
Thus if the limit~\eqref{equ:Lyapunov_exp} exists, then for large $n$ we have that
\[ 
e^{\lambda(x_{0})} \approx \sqrt[n]{ \prod_{i=0}^{n}\|h'(x_{i})\| }.
\]
So $e^{\lambda(x_{0})}$ is an average value of the norm of the tangent map along the orbit of $x_0$.

In particular, if $\lambda(x_0)<0$, then (saying non-strictly) \myemph{``in average the orbit of $x_0$ should attract points close to it''}.
Hence the points with negative Lyapunov exponents would detect attractors of dynamical systems.

Moreover, suppose there is a diffeomorphism $h:\Mman\to\Mman$ isotopic to the identity (e.g. a diffeomorphism belonging to a flow) which also have negative Lyapunov exponents at some large subsets of $\Mman$.
Then one would expect that representatives of certain (co)homology classes of $\Mman$ can be deformed under iterations of $h$ into subsets of small sizes, and therefore those classes could vanish.
In other words, one would get triviality of some (co)homology or homotopy groups of $\Mman$.

Stochastic analogues of Lyapunov exponents, the so called \myemph{$p$-moment exponents}, play a crucial role for investigation of stability of stochastic dynamical systems, see e.g. \cite{Khasminskii:JAMM:1962}, \cite{Kozin:JMAA:1965}, \cite{KozinSugimoto:Proc:1977}, \cite{KozubovskajHrisanov:UMJ:1981},  \cite{LArnold:1984}, \cite{ArnoldKliemannOeljeklaus:LNM:1986}, \cite{BaxendaleStroock:PTRF:1988}, \cite{ElwRos:1996-Hom} and others.

Let $\xsol$ be a stochastic flow on $\Mman$ being a solution of SDE with smooth coefficients.
Roughly speaking it is a family of differentiable flows depending on a some parameter $\omega$ belonging to a probability space $\Omega$, see Definition~\ref{def:cont_stoch_flow} below.
Then given a compact subset $K\subset \Mman$ and $p>0$ define the \myemph{$p$-th moment exponent} by
\[
\mu_{K}(p) := 
\overline{\lim_{t\to\infty}} \ \sup_{x\in K} \frac{1}{t} \ln \EE \| T_x \xsol_t \|^p.
\]
A stochastic flow $\xi$ is called \myemph{$p$-th moment stable} whenever $\mu_{x}(p)<0$ for each $x\in\Mman$ and \myemph{strongly $p$-th moment stable} if $\mu_{K}(p)<0$ for each compact subset $K\subset\Mman$.

Then the strong $p$-th moment stability would imply that in average the flow decreases the sizes of compact sets.
In particular, it was shown in Elworthy and Rosenberg~\cite{ElwRos:1996-Hom} that for a compact manifold $\Mman$
\begin{itemize}
\item
$\mu_{\Mman}(1)<0$ implies triviality of the fundamental group of $\Mman$ vanish: $\pi_1\Mman=0$;
\item
$\mu_{\Mman}(2)<0$ implies triviality of the second homotopy group of $\Mman$ vanish: $\pi_2\Mman=0$;
\item
$\mu_{\Mman}(q)<0$ implies triviality of $q$-th homologies of $\Mman$: $H_q(\Mman,\ZZZ)=0$;
\item
if $\mu_{\Mman}([\tfrac{n+1}{2}])<0$, then $\Mman$ is a homotopy sphere.
\end{itemize}

For non-compact manifolds the situation is more complicated as a priori one can not expect uniform bounds for $\mu_K(p)$.
That case was considered by Xue Mei Li~\cite{XMLi:PTRF:1994}.
She studied moment stability of SDE of the form
\begin{equation}\label{equ:SDE_Li}
dx_t = X(x_t) \circ dB_t + A(x_t) dt,
\end{equation}
where $B_t$ is an $m$-dimensional Brownian motion on $\mathcal{T}$, $A$ is a vector field on $\Mman$, and $\Xmap\in Hom(\Rntriv, T\Mman)$ is a bundle homomorphism of class $C^3$ from trivial $\RRR^n$-bundle $\Rntriv=\RRR^n\times\Mman\to\Mman$ over $\Mman$ to its tangent bundle $T\Mman\to\Mman$.
Among other results she gave sufficient conditions for triviality of the fundamental group $\pi_1\Mman$ of the manifold $\Mman$ in terms of the coefficients of~\eqref{equ:SDE_Li}.

The aim of the present paper is to show that the method of proof of Theorem~4.1 in~\cite{XMLi:PTRF:1994} allows to establish an essentially stronger result about topological structure of $\Mman$: namely that all the homotopy groups $\pi_n\Mman$ vanish, so $\Mman$ is \myemph{contractible}, see Theorem~\ref{th:M_is_contractible1} below.

\section{Preliminaries}
We start with a usual setting of theory of SDE.
Let $\Mman$ be a smooth connected manifold (i.e.\! locally Euclidean Hausdorff topological space with countable base) of dimension $m$ possibly non-compact and having a boundary and $\mathcal{T}=(\Omega, \FF, \PPP)$ be a probability space, so $\Omega$ is a set, $\FF$ is a $\sigma$-algebra of subsets of $\Omega$, and $\PPP$ is a probability measure on $\FF$.
Let also $\{\FF_t\}_{t\geq0}$ for some $a\geq0$ be a family of $\sigma$-algebras in $\FF$ with the following properties:
\begin{enumerate}[label={\rm(\roman*)}]
\item\label{enum:sigma_algebras:null_sets}
each $\FF_t$ contains all null sets of $\FF$;
\item\label{enum:sigma_algebras:increasing}
$\FF_s \subseteq \FF_t$ for $s<t$;
\item\label{enum:sigma_algebras:right_continuous}
$\{\FF_t\}_{t\geq0}$ is right continuous in the sense that $\FF_s = \cap_{s<t} \FF_t$ for all $s\geq0$.	
\end{enumerate}	
For a topological space $Y$ we will denote by $\mathcal{B}(Y)$ the Borel $\sigma$-algebra of subsets of $Y$.
Given a map $f:A\times B \times C \to D$ of a product of topological spaces we will often consider \myemph{restriction} maps obtained by fixing some coordinates, e.g. $f_{a}:\{a\}\times B\times C\to D$ defined by $f_a(b,c) = f(a,b,c)$, or $f_{a,b}:\{a\}\times \{b\} \times C\to D$, $f_{a,b}(c) = f(a,b,c)$, for $(a,b,c)\in A\times B \times C$.
Thus we put the corresponding fixed coordinates as subindexes.

Denote $\Delta = \{(s,t)\in\RRR^2\}\mid 0 \leq s \leq t \}$.

\begin{definition}\label{def:cont_stoch_flow}
By a \myemph{continuous stochastic flow} we mean a map
\begin{equation}\label{equ:cont_stoch_flow}
\xsol:\Mman\times\Delta\times\Omega \to \Mman
\end{equation}
having the following properties: there exists a subset $N\in\mathcal{F}$ of measure $0$ such that for all $x\in\Mman$, $(s,t)\in\Delta$, and $\omega\in\Omega\setminus N$
\begin{enumerate}[label={\rm\alph*)}]
\item\label{enum:def:cont_stoch_flow:xi_xt}
the map $\xsol_{x,s,t}:\Omega\to\Mman$, $\xsol_{x,s,t}(\omega) = \xsol(x,s,t,\omega)$, is $\FF_t/\mathcal{B}(\Mman)$-measurable;
\item\label{enum:def:cont_stoch_flow:xi_omega}
the map $\xsol_{\omega}:\Mman\times\Delta\to\Mman$, $\xsol_{\omega}(x,s,t) = \xsol(x,s,t,\omega)$, is continuous;
\item\label{enum:def:cont_stoch_flow:xi_ss_id}
$\xsol_{s,s,\omega}(x)=x$;
\item\label{enum:def:cont_stoch_flow:xi_st__xi_tu}
$\xsol_{t,u,\omega}(\xsol_{s,t,\omega}(x)) = \xsol_{s,u,\omega}(x)$, whenever $0\leq s\leq t \leq u$,
\end{enumerate}
where the map $\xsol_{t,u,\omega}:\Mman\to\Mman$ is given by $\xsol_{t,u,\omega}(x) = \xsol(x, t,u,\omega)$.
\end{definition}

\begin{definition}\label{def:stoch_deformation}
By a \myemph{stochastic deformation} we will mean a map
\begin{equation}\label{equ:stoch_deformation}
\axsol:\Mman\times [0,+\infty)\times\Omega \to \Mman
\end{equation}
having the following properties: there exists $N\in\FF$ of measure $0$ such that for each $\omega\in\Omega\setminus N$
\begin{enumerate}[label={\rm\alph*)}]
\item\label{enum:def:stoch_deformation:xi_xt}
the map $\xsol_{x,t}:\Omega\to\Mman$, $\xsol_{x,t}(\omega) = \xsol(x,t,\omega)$, is $\FF_t/\mathcal{B}(\Mman)$-measurable;
\item\label{enum:def:stoch_deformation:xi_omega}
the map $\axsol_{\omega}:\Mman\times [0,+\infty) \to\Mman$, $\axsol_{\omega}(x,t) = \axsol(x,t,\omega)$, is continuous;
\item\label{enum:def:stoch_deformation:xi_s_id}
$\axsol(x,0,\omega) = x$ for all $x\in\Mman$.
\end{enumerate}
If, in addition, the following condition holds:
\begin{enumerate}[label={\rm\alph*)}, resume] 
\item\label{enum:def:stoch_deformation:xi_st__xi_tu}
$\axsol_{t,\omega}(\axsol_{s,\omega}(x)) = \axsol_{t+s,\omega}(x)$, for all $s,t\geq0$
\end{enumerate}
then $\axsol$ will be called an \myemph{autonomous continuous stochastic flow}.
\end{definition}

\begin{remark}\rm
1) Every continuous stochastic flow $\xsol$ defines a family $\{\axsol_s\}_{s\geq0}$ of stochastic deformations 
\begin{align*}
&\axsol_s:\Mman\times [0,+\infty)\times\Omega \to \Mman, &
\axsol_s(x,t,\omega) &= \xsol(x, s, s+t, \omega).
\end{align*}

2) Every autonomous continuous stochastic flow $\axsol$ defines a continuous stochastic flow 
\begin{align*}
&\xsol:\Mman\times\Delta\times\Omega \to \Mman, &
\xsol(x,s,t,\omega) &= \axsol(x,t-s,\omega)
\end{align*}
satisfying $\xsol_{s,t,\omega}=\xsol_{s+\tau,t+\tau,\omega}$ for all $\tau\geq0$.

3) Every stochastic deformation $\axsol$ such that $\axsol_{t,\omega}:\Mman\to\Mman$ is a homeomorphism for all $t\geq0$ and $\omega\in\Omega\setminus N$ defines a continuous stochastic flow 
\begin{align}\label{equ:stoch_def_hom__stoch_flow}
&\xsol:\Mman\times\Delta\times\Omega \to \Mman, &
\xsol(x,s,t,\omega) &= \axsol_{t,\omega}^{-1}\circ\axsol_{s,\omega}(x).
\end{align}
Let us mention that for a stochastic deformation we do not require the map $\axsol_{t,\omega}$ to be invertible.
\end{remark}

It is well known that for a large class of (autonomous) SDE their solutions are (autonomous) continuous stochastic flows in the above sense, however a priori, not every stochastic flow is a solution of certain SDE.
For details on this correspondence, see e.g. \cite[Chapter~4]{Kunita:SF_SDE:1990}, \cite{Kunita:St:1981},  \cite{Taniguchi:SSR:1989}.

\subsection{Homotopies}
We will briefly recall the necessary definitions.
Let $S$ and $\Mman$ be two topological spaces and $f,g:S \to \Mman$ be two continuous maps between them.
These maps are called \myemph{homotopic} if there exists a (jointly) continuous map $H:S \times I \to \Mman$ such that $H(0,x) = f(x)$ and $H(x,1) = g(x)$ for all $x\in S$.
Any such map $H$ is called a \myemph{homotopy} between $f$ and $g$.

A map $f:S \to \Mman$ homotopic to a constant map is also said to be \myemph{null homotopic}.

A \myemph{deformation} of $\Mman$ is a homotopy $H:M\times I \to M$ starting from $\id_{\Mman}$, i.e. $H_0 = \id_{\Mman}$.
Thus, roughly speaking, a \myemph{stochastic deformation} in the sense of Definition~\ref{def:cont_stoch_flow} is a family of deformations ``measurably depending'' on some parameter $\omega\in\Omega$.

A topological space $\Mman$ is said to be \myemph{contractible} if the identity map $\id_{\Mman}:\Mman\to\Mman$ is null-homotopic, i.e. homotopic to a constant map $*:\Mman\to p \in \Mman$ into some point $p\in \Mman$.
The corresponding homotopy between $\id_{\Mman}$ and $*$, i.e.\! a ``deformation of $\Mman$ into a point $p$'', is called a \myemph{contraction of $\Mman$}.

For instance, \myemph{any convex subset $\Mman \subset \RRR^n$ is contractible}.
Indeed, let $p\in\Mman$ be any point, then a contraction $H:\Mman\times I \to \Mman$ of $\Mman$ into $p$ can be given by the formula: $H(x,t) = t p  + (1-t)x$.
On the other hand each compact manifold without boundary, e.g. the $n$-dimensional sphere $S^n$ and the $n$-torus $T^n$ are not contractible.

Notice that \myemph{if $\Mman$ is contractible, then each continuous map $\sigma: S \to \Mman$ is null homotopic}.
Indeed, if $H:\Mman\times I \to \Mman$ is a contraction of $\Mman$ into some point $p \in \Mman$, then the map $\Sigma:S \times I \to \Mman$ defined by $\Sigma(x,t) = H(\sigma(x), t)$ is a homotopy between $\sigma$ and a constant map into the point $p$.

The following statement is a particular case of the well-known Whitehead's theorem.

\begin{theorem}\label{th:contractibility_criterium}{\rm (J.~H.~C.~Whitehead), e.g.~\cite[Theorem~4.5]{Hatcher:AT:2002}.}
A connected manifold $\Mman$ is contractible if and only if for each $n\geq1$ each continuous map $\sigma:S^n\to \Mman$ is null homotopic.
\end{theorem}

\subsection{Measures associated with a stochastic deformation.}
Let $\xsol$ be a stochastic deformation on $\Mman$.
By assumption~\ref{enum:def:stoch_deformation:xi_xt} of Definition~\ref{def:stoch_deformation}, for each $(x,t)\in\Mman\times[0,+\infty)$ the map $\xsol_{x,t}$
\begin{align*}
&\xsol_{x,t}:\ \Omega \ \to \ \Mman, &
\xsol_{x,t}(\omega) = \xsol(x,t,\omega).
\end{align*}
is $\FF_t/\mathcal{B}(\Mman)$-measurable.
Therefore one can define the following $\sigma$-additive probability measure $\mu_{x,t}$ on $\Mman$ by
\[
\mu_{x,t}(K) := \PPP(\xsol_{x,t}^{-1}(K)) = \PPP\{ \omega\in\Omega \mymid \xsol_{t,\omega}(x)\in K \},
\qquad K \in \mathcal{B}(\Mman).
\]

Our aim is to prove the following theorem.
\begin{theorem}\label{th:M_is_contractible1}
Let $\Mman$ be a smooth connected complete Riemannian manifold and $\mathcal{T}=(\Omega, \FF, \PPP)$ be a probability space.
Suppose there exists a stochastic deformation
\[
\xsol:\Mman\times [0,+\infty)\times\Omega \to \Mman
\]
having the following properties:
\begin{enumerate}[leftmargin=*, label=\rm(\roman*)]
\item\label{enum:cond:xi_is_diffeo}
the map $\xsol_{t,\omega}$ is $C^1$ for all $t\in[0,+\infty)$ and $\omega\in\Omega\setminus N$, where $N$ is a null set from Definition~\ref{def:stoch_deformation};
\item\label{enum:cond:Lie_condition}
for each compact subset $K\subset\Mman$  we have that
\begin{equation}\label{equ:Li_condition1}
\int\limits_{0}^{+\infty}\,
\sup\limits_{x\in K}\,\EE\,\|T_x \xsol_{t,\omega}\| dt < \infty,
\end{equation}
where the norm is taken with respect to the corresponding Riemannian metric;

\item\label{enum:cond:muK}
there exist a point $z\in\Mman$ and a compact subset $K\subset\Mman$ such that 
\[
\varliminf_{t\to\infty} \mu_{z,t}(K) \equiv
\varliminf_{t\to\infty} \PPP\{ \omega\in\Omega \mymid \xsol_{t,\omega}(z)\in K \}
 > 0,
\]
i.e. one can find $\varepsilon>0$ and $A>0$ satisfying $\mu_{z,t}(K)>\varepsilon$ for all $t>A$.
\end{enumerate}
Then $\Mman$ is \myemph{contractible}.
\end{theorem}

\section{Proof of Theorem~\ref{th:M_is_contractible1}}\label{sec3}

\subsection{}
Recall that a complete Riemannian manifold $\Mman$ is also a complete metric space with the distance $\rho(x,y)$ between points $x,y\in\Mman$ defined as the infimum of lengths of $C^1$-paths $\gamma:[0,1]\to\Mman$ such that $\gamma(0)=x$ and $\gamma(1)=y$.

Also for each $x\in\Mman$ and a unit tangent vector $v\in T_x\Mman$ there exists a unique geodesic line $\gamma_{x,v}:\RRR\to\Mman$ such that $\gamma_{x,v}(0)=x$ and $\dot{\gamma}_{x,v}(0)=v$.
This allows to define the following \myemph{exponential} map $\exp_x: T_x\Mman\to\Mman$ by
\[
\exp_x(w) = \gamma_{x,\frac{w}{|w|}}(|w|), 
\qquad 
w\in T_x \Mman,
\]
so $\exp_x(x)=0$ and it maps radial lines $\{tv\}_{t\in\RRR}$ onto geodesics passing through $x$.

It is well known that $\exp_x$ is a $C^{\infty}$ map being a local diffeomorphism at $0\in T_x\Mman$.
Let $D_r(0) \subset T_x\Mman$ be an open ball of radius $r$ with center at the origin.
Then there exists $r>0$ such that $\exp_x$ diffeomorphically maps $D_r(0)$ onto some neighbourhood $B_r(x)$ of $x$.
Such a neighborhood $B_r(x)$ is called a \myemph{geodesic ball} at $x$ of radius $r$ and the supremum of all such $r$ for which $B_r(x)$ is defined is called the \myemph{injectivity} radius at $x$ with respect to $\rho$ and denoted by $\Rinj{x}$.
Thus
\[ \Rinj{x} = \sup_{r>0} \{r \mymid  \text{the restriction} \ \exp_x|_{D_r(0)}:D_r(0) \to \Mman \ \text{is an embedding}\}.\]
For a subset $K\subset\Mman$ let 
\begin{equation}\label{equ:inj_radius}
\Rinj{K} = \inf_{x\in K} \Rinj{x}.
\end{equation}
If $K$ is compact, then $\Rinj{K}>0$.

\subsection{}
The proof of contractibility of $\Mman$ follows the line of~\cite[Theorem~4.1]{XMLi:PTRF:1994}.
Due to Whitehead Theorem~\ref{th:contractibility_criterium} it suffices to show that every continuous map $\sigma:S^n \to \Mman$ from $n$-dimensional sphere $S^n$ into $\Mman$ is null homotopic, i.e.\! homotopic to a constant map.

Since every continuous map $\sigma:S^n\to\Mman$ is homotopic to a $C^{1}$-map, we may assume that \myemph{$\sigma$ is of class $C^{1}$}.

So let $\sigma:S^n\to M$, $(n\geq1)$, be a $C^1$ map.
For each $\omega\in\Omega\setminus N$ and $t\in[0,+\infty)$ define the map $\sigma_{t,\omega}: S^n\to\Mman$ by
\[
\sigma_{t,\omega}(x) = \xsol_{t,\omega}(\sigma(x)).
\]
Then due to~\ref{enum:cond:xi_is_diffeo} \myemph{for all $\omega\in\Omega\setminus N$ and $t\in[0,+\infty)$ the map $\sigma_{t,\omega}=\xsol_{t,\omega}\circ\sigma$ is $C^1$} as well.

\begin{lemma}\label{lm:reduction_of_theorem}
Consider the following conditions:
\begin{enumerate}[leftmargin=*, label=$(\arabic*)$]
\item\label{enum:lm:reduction_of_theorem:sigma_null_homotopic}
$\sigma$ is null homotopic;
\item\label{enum:lm:reduction_of_theorem:sigma_tomega_null_homotopic}
there exist $\omega\in\Omega\setminus N$ and $t>0$ such that $\sigma_{t,\omega}$ is null homotopic;
\item\label{enum:lm:reduction_of_theorem:sigma_in_ball}
 there exist $\omega\in\Omega\setminus N$ and $t>0$ such that $\sigma_{t,\omega}$ is contained in a geodesic ball $B$ near some point $x\in\Mman$;
\item\label{enum:lm:reduction_of_theorem:set_Z}
there exists a subset $Z\subset\Omega\setminus N$ of positive measure such that for each $\omega\in Z$ there exists $t>0$ such that $\sigma_{t,\omega}$ is contained in a geodesic ball $B$ near some point $x\in\Mman$.
\end{enumerate}
Then we have the following implications: 
\ref{enum:lm:reduction_of_theorem:set_Z} $\Rightarrow$
\ref{enum:lm:reduction_of_theorem:sigma_in_ball} $\Rightarrow$
\ref{enum:lm:reduction_of_theorem:sigma_tomega_null_homotopic} $\Rightarrow$
\ref{enum:lm:reduction_of_theorem:sigma_null_homotopic}.
\end{lemma}
\begin{proof}
Suppose~\ref{enum:lm:reduction_of_theorem:set_Z} holds.
Since $Z$ has positive measure, it is non-empty, which implies~\ref{enum:lm:reduction_of_theorem:sigma_in_ball}.

Suppose~\ref{enum:lm:reduction_of_theorem:sigma_in_ball} holds.
As each geodesic ball $B_r$ is homeomorphic (e.g. via the exponential map) with a standard $n$-dimensional ball $D_r$ and $D_r$ is convex and therefore contractible, we obtain that $\sigma_{t,\omega}:S^n \to B_r \subset \Mman$ is null-homotopic, that is~\ref{enum:lm:reduction_of_theorem:sigma_tomega_null_homotopic} is satisfied.

Finally, suppose~\ref{enum:lm:reduction_of_theorem:sigma_tomega_null_homotopic} holds.
Since $\xsol_{\omega}:\Mman$ is homotopic to the identity map $\xsol_{0,\omega} = \id_{\Mman}$, whence $\sigma_{t,\omega}=\xsol_{t,\omega}\circ\sigma$ is homotopic to $\sigma=\xsol_{0,\omega}\circ\sigma$.
Therefore if $\sigma_{t,\omega}$ is null homotopic, then so is $\sigma$.
\end{proof}

Thus for the proof of Theorem~\ref{th:M_is_contractible1} it suffices to find a set $Z$ satisfying~\ref{enum:lm:reduction_of_theorem:set_Z} of Lemma~\ref{lm:reduction_of_theorem}.

\begin{lemma}\label{lm:lim_E_diam_st_0}
For each $t\in[0,\infty)$ let $\diam(\sigma_{t}):\Omega\to[0,+\infty)$ be the random variable equals the diameter of the image of $\sigma_{t,\omega}(S^n)$ in $\Mman$ with respect to the metric $\rho$.
Then there exists a sequence of numbers $\{t_j\}\subset\RRR$ converging to infinity and such that 
\begin{equation}\label{equ:lim_E_diam_st_0}
\lim_{j\to\infty} \EE \bigl(\diam(\sigma_{t_j})\bigr)  = 0.
\end{equation}
\end{lemma}
Assuming that Lemma~\ref{lm:lim_E_diam_st_0} is proved we will complete Theorem~\ref{th:M_is_contractible1}.

Let $z$ and $K$ be the same as in~\ref{enum:cond:muK} of Theorem~\ref{th:M_is_contractible1}, so there exists $\varepsilon>0$ and $A>0$ such that 
\begin{equation}\label{equ:P_greater_muK2}
\mu_{x,t}(K) \ = \ \PPP\{ \omega\in\Omega \mymid \xsol_{t,\omega}(z) \in K\} \ > \ \varepsilon
\end{equation}
for all $t>A$.

Denote by $\Rinj{K}$ the injectivity radius of $K$ with respect to the metric $\rho$.
Then $\Rinj{K}>0$.

Let also $\{t_j\}\subset(0,+\infty)$ be a sequence satisfying~\eqref{equ:lim_E_diam_st_0}.
Then there exists $t_j > A$ such that
\begin{equation}\label{equ:P_less_muK4}
\PPP \{\omega\in\Omega\mymid \diam(\sigma_{t_{j}})(\omega) \geq \Rinj{K}/2 \} \ < \ \varepsilon/2.
\end{equation}

Consider the following set  
\[
Z = \{\omega\in\Omega \mymid \diam(\sigma_{t_j})(\omega) < \Rinj{K}/2 \ \text{and} \ \xsol_{t_j,\omega}(z) \in K  \}.
\]
Then for every $\omega\in Z$ we have that the image of $\sigma_{t_j,\omega}(S^n)=\xsol_{t_j,\omega}\circ\sigma(S^n)$ intersects $K$ and is contained in the geodesic ball of radius smaller than $\Rinj{K}$ with center $\xsol_{t_j,\omega}(z)$.
Therefore the map $\sigma_{t_j,\omega}:S^n\to K \subset \Mman$ is null homotopic, whence so is $\sigma$.

It remains to show that $\PPP(Z)>0$.
Indeed,
\begin{align*}
\PPP(Z) &= \PPP \{ \omega\in\Omega \mymid  \xsol_{t_j}(z) \in K \}  - 
       \PPP \{ \omega\in\Omega \mymid \diam(\sigma_{t_j}(\omega) \geq \Rinj{K}/2 \ \text{and}\ \xsol_{t_j}(z) \in K \} \\
&\geq \PPP \{ \omega\in\Omega \mymid \xsol_{t_j}(z) \in K \} - \PPP\{\omega\in\Omega \mymid \diam(\sigma_{t_j})(\omega) \geq \Rinj{K}/2\} \\
&> \varepsilon \ - \ \varepsilon/2 \ = \ \varepsilon/2 \ > \ 0.
\end{align*}
The third line is obtained due to~\eqref{equ:P_greater_muK2} and~\eqref{equ:P_less_muK4}.
Thus $Z$ is non-empty and satisfies~\ref{enum:lm:reduction_of_theorem:set_Z} of Lemma~\ref{lm:reduction_of_theorem}.
This proves that $\Mman$ is contractible modulo Lemma~\ref{lm:lim_E_diam_st_0}.

\subsection{Proof of Lemma~\ref{lm:lim_E_diam_st_0}}
Notice that every great circle $e$ in $S^n$ of radius $1$ can be regarded as a length preserving map $e:[0,2\pi]\to S^n$ and we will denote by $\dot{e} = \partial e / \partial s$ the unit tangent vector field along $e$.

This circle is uniquely determined by a $2$-plane in $\RRR^{n+1}$ passing through the origin and so the space of all great circles in $S^n$ can be identified with the Grassmannian manifold $G^{n+1}_{2}$ of $2$-planes in $\RRR^{n+1}$.

For each $t\in[0,+\infty)$ and a great circle $e\in G^{n+1}_{2}$ define the map $l_{t,e}:\Omega\setminus\Nman\to[0,+\infty)$ associating to each $\omega\in\Omega\setminus\Nman$ the length of the  curve $\sigma_{t,\omega}\circ e$.
Thus
\begin{equation}\label{equ:length_lte}
\begin{aligned}
l_{t,e}(\omega) &= length (\sigma_{t,\omega}\circ e) \\
&= \int\limits_{0}^{2\pi} \Bigl|T_{e(s)}\bigl(\xsol_{t,\omega} \circ \sigma\bigr)(\dot{e}(s))\Bigr| ds = 
\int\limits_{0}^{2\pi} \Bigl|T_{\sigma(e(s))} \xsol_{t,\omega} \circ T_{e(s)}\sigma (\dot{e}(s))\Bigr| ds.
\end{aligned}
\end{equation}

Notice that if $x,y\in S^n$ and $e$ be a great circle passing through $x$ and $y$, then for each $\omega\in\Omega\setminus\Nman$ we have that 
\[
\rho(\sigma_{t,\omega}(x), \sigma_{t,\omega}(y)) \leq \tfrac{1}{2}\, l_{t,e}(\omega),
\]
whence 
\begin{equation}\label{equ:diam_less_than_2_length}
\diam(\sigma_{t})(\omega) \ \leq \ \tfrac{1}{2}\, \sup_{\text{$e$ is a great circle}} l_{t,e}(\omega).
\end{equation}

Let $p:TS^n\to S^n$ be the tangent bundle of $S^n$ and $US^n \subset TS^n$ be the sphere-bundle consisting of all tangent vectors of length $1$.
Denote by $\mathbf{L} = T\sigma(US^n) \subset TM$ the image of $US^n$ in $TM$ under the tangent map $T\sigma:TS^n\to TM$.
Evidently, $US^n$ and $\mathbf{L}$ are compact.

Notice that if $e:[0,2\pi]\to S^n$ is a great circle in $S^n$ and $\dot{e} = \partial e / \partial s$, then
\begin{align*}
\bigl(e(s), \dot{e}(s)\bigr) &\in US^n, &
\bigl(\, \sigma\circ e(s), \ T_{e(s)}\sigma (\dot{e}(s)) \, \bigr) & \in  \mathbf{L}.
\end{align*}
Therefore we get from~\eqref{equ:Li_condition1} and~\eqref{equ:length_lte} that
\begin{align*}
\int\limits_{0}^{+\infty}\, \sup_{e \in G^{n+1}_{2}} \, \EE \, l_{t,e} dt 
&= \int\limits_{0}^{+\infty}
\sup_{e \in G^{n+1}_{2}}\, \EE \,  
\int\limits_{0}^{2\pi}\Bigl|T_{\sigma(e(s))} \xsol_{t,\omega} \circ T_{e(s)}\sigma (\dot{e}(s))\Bigr| ds\, dt \\
&= \int\limits_{0}^{+\infty}
\sup_{e \in G^{n+1}_{2}}\, 
\int\limits_{0}^{2\pi} \, \EE \,  \Bigl|T_{\sigma(e(s))} \xsol_{t,\omega} \circ T_{e(s)}\sigma (\dot{e}(s))\Bigr| ds\, dt \\
&= 2\pi \int\limits_{0}^{+\infty}
\sup_{e \in G^{n+1}_{2}}\, 
\, \sup_{s\in[0,2\pi]} 
\, \EE \,  \Bigl|T_{\sigma(e(s))} \xsol_{t,\omega} \circ T_{e(s)}\sigma (\dot{e}(s))\Bigr|  dt \\
&\leq 2\pi \int\limits_{0}^{+\infty} 
  \sup_{e \in G^{n+1}_{2}} \,
  \sup_{s\in[0,2\pi]} \,
  \EE \,
  \Bigl( \| T_{\sigma(e(s))} \xsol_{t,\omega} \|
  \sup_{v\in\mathbf{L}\, \cap\, T_{\sigma(e(s))}\Mman }\, |v| \Bigr) \, dt  \\
& \leq  
2\pi C \int\limits_{0}^{+\infty} \, \sup_{x\in \sigma(S^n)} \, \EE \, \| T_{x} \xsol_{t,\omega} \|  dt
\stackrel{\eqref{equ:Li_condition1}}{<} \infty,
\end{align*}
where
\[
C = \sup_{(x,v)\in \mathbf{L}} |v| 
\]
is finite due to compactness of $\mathbf{L}$.
It follows that $\lim\limits_{t\to\infty} \sup\limits_{e \in G^{n+1}_{2}} \, \EE \, l_{t,e} = 0$.
Therefore we get from~\eqref{equ:diam_less_than_2_length} that there exists a sequence $\{t_j\}_{i\geq0}$ such that $\lim\limits_{j\to\infty} \EE \bigl(\diam(\sigma_{t_j})\bigr)  = 0$.

This completes proof of Lemma~\ref{lm:lim_E_diam_st_0} and therefore Theorem~\ref{th:M_is_contractible1} as well.
\hfill\qed


\end{document}